\documentstyle{article}

\newtheorem{theorem}{Theorem}
\newtheorem{lemma}[theorem]{Lemma}
\newtheorem{proposition}[theorem]{Proposition}

\def\eps{\varepsilon}
\def\a{\alpha}

\def\s{\sigma}
\def\g{\gamma}
\def\k{\kappa}
\def\G{\Gamma}
\def\E{{\bf E}}
\def\P{{\bf P}}
\def\Q{{\bf Q}}
\def\k{\kappa}
\def\tG{\tilde{\G}}

\title{Average number of distinct  part sizes   in a random Carlitz
composition} 

\author{William M. Y. Goh
\\
{\small  Department of Mathematics and Computer Science}\\
{\small Drexel University}\\
{\small Philadelphia, PA 19104}\\
{\small \texttt{wgoh@mcs.drexel.edu}}\\
\and 
Pawe{\l} Hitczenko       \\
{\small  Department of Mathematics and Computer Science}\\
{\small Drexel University}\\
{\small Philadelphia, PA 19104}\\
{\small \texttt{phitczen@mcs.drexel.edu}}\\
}

\begin{document}
 
\begin{titlepage} 
\date{ \today}
 
\maketitle
 
\abstract{A composition of an integer $n$ is called Carlitz if 
adjacent parts are different. Several characteristics of random
Carlitz compositions have been studied recently by Knopfmacher and
Prodinger. We will complement their work by establishing 
asymptotics of the average
number of distinct part sizes in a random Carlitz composition.}
\end{titlepage}

\section{Introduction}

In this  note we obtain  precise  asymptotics, as $n\to\infty$, for the
expected number of distinct 
 part sizes in a random Carlitz composition  of an integer $n$. Let us recall  that
a tuple $(\g_1,\dots,\g_k)$ is a {\em composition} of an integer $n$ 
if the $\g_j$'s are positive integers,  called parts,  such that
$\sum_j\g_j=n$. The number $k$ is the number of parts and the
values of 
$\g_j$'s are called part sizes. There are $2^{n-1}$different
compositions of $n$. A composition is called {\em Carlitz}
if the adjacent parts are different, i.e. if $\g_j\ne\g_{j+1}$ for $j=1,\dots,k-1$.
 For example, out of  sixteen compositions of the integer 5, seven
are 
Carlitz, namely $(5)$, $(4,1)$, $(1,4)$, 
$(3,2)$, $(2,3)$, $(1,3,1)$,  and $(2,1,2)$. For other values, see \cite[Sequence A003242]{S}. We denote the set of all
Carlitz compositions of $n$ by $\Omega_n$. Carlitz compositions have
been introduced by Carlitz \cite{Carlitz} who found the generating function for the
total number of them. They have been subsequently studied by
Knopfmacher and Prodinger \cite{KP} (see also \cite{LP}). These authors found the
asymptotics of the total number of Carlitz compositions. They also
studied several parameters (like the number of parts, the size of the
largest part among other things) for {\em random} Carlitz
compositions. \lq\lq Random Carlitz composition\rq\rq\  means  a
composition  chosen according to the  uniform probability measure on
$\Omega_n$. This measure will be denoted by $\P$ and $\E$ will denote
integration with respect to $\P$. In
this setting, various parameters of Carlitz compositions become
 random variables and their probabilistic properties are to be
studied. For example Knopfmacher and Prodinger found, among other things, the exact
asymptotic behavior of the expected number of
parts and the expected size of the largest part. One question that
they left open concerned the expected value of the number of {\em distinct}
part sizes, $D_n$, and the purpose of this note is to answer their
question.
 In order to state the result we need to introduce some more
notation and we will try to closely follow the notation of Knopfmacher and
Prodinger.  First, the  number of 
distinct part sizes is defined
formally as follows: if $(\g_1,\dots,\g_k)$ is in
$\Omega_n$ then
$$D_n=1+\sum_{i=2}^kI_{\{\g_i\ne\g_j,\ j=1,\dots,i-1\}},$$
where $I_A$ denotes, as usual, the indicator function of the set $A$.   
Secondly, we introduce a function $\sigma$  of a complex variable
defined by 
$$\sigma(z)=\sum_{j=1}^\infty(-1)^{j-1}\frac{z^j}{1-z^j}.$$
The equation $\sigma(z)=1$ has the unique real solution $\rho=0.571349\dots$
 on the
interval $[0,1]$.  (The relevance of
this is that, as was shown by Carlitz, the generating function of
Carlitz compositions is equal to $1/(1-\sigma(z))$.) We have
\begin{theorem} With the above notation, and letting $\{\ \cdot\ \}$ denote the fractional part, as $n\to\infty$, we have 
$$\E D_n=\frac{\ln( n/\s'(\rho))}{\ln(1/\rho)}+\frac12-\frac{\g}{\ln(1/\rho)}+h_0(\{\rho\ln(n/\s'(\rho)\})+o(1),$$
where $\g$ is Euler's constant, $h_0$ is a mean zero function of
 period 1 whose Fourier coefficients are given by
$$
c_\ell=
\frac1{\ln(1/\rho)}\G(-\frac{2\pi i\ell}{\ln(1/\rho)}),\quad \ell\ne0.
$$
\end{theorem}

Approximating all the constants  and using the fact that the gamma function decays very fast along the imaginary axis, it can
be seen in particular, 
that 
$$\E D_n=C_1\ln n-C_2+h_0(\{\rho\ln(n/\s'(\rho)\})+o(1),$$
where $C_1=1.786495\dots$, $C_2=2.932545\dots$, and the amplitude of $h_0$  is bounded by $0.5882304\cdots\times 10^{-7}$.

The approach is as in \cite{KP} via generating functions. We let $I_j$ denote the set of those Carlitz
compositions that  contain at least one  part of size \lq\lq
j\rq\rq.
Using, without any risk of confusion, the same notation for a
set and its indicator we have
$$D_n=\sum_{j=1}^nI_j,$$ and therefore
\begin{equation}
\E D_n=\sum_{j=1}^n\P(I_j)=\sum_{j=1}^n(1-\P(I_j^c)),\label{edn}
\end{equation}
where $A^c$ denotes the complement of a set $A$.
Denoting by $a_n$ and $a_{n,j}$ the number of all Carlitz compositions of $n$
and the number of those Carlitz compositions of $n$ that {\em do not}
use size $j$ we have 
$$\P(I_j^c)=\frac{a_{n,j}}{a_n}.$$ The sequence $(a_n)$ was studied in
Knopfmacher and Prodinger \cite{KP}, so we need only to study the numbers
$a_{n,j}$. In order to do that we will build their generating
function. The construction follows the ideas of \cite{KP}, but since that paper is short
on some details and may be a bit difficult to read for a new adept, we will provide a fairly detailed argument.

We would like to mention that the method of our paper does not appear to be strong enough to yield  information about the limiting distribution of the variable $D_n$.  A bivariate generating function for $D_n$ would be extremely welcome. But, we were unable to get it. Although it is known (c.f. \cite{LP}) that the total number of parts in Carlitz composition satisfies the central limit theorem, we think that it is unlikely that $D_n$ will have the same property. It would be very interesting and desirable to find the limiting distribution of $D_n$.

\section{Generating function}
In this section we will prove the following statement 
\begin{proposition} Let $C_j(z)$ be the generating function of the sequence $\{a_{n,j},\ n\ge0\}$. Then,
$$C_j(z)=\frac1{1-\s_j(z)},$$
where
\begin{equation}\s_j(z)=\sum_{\ell\ge1}(-1)^{\ell-1}\left(\frac{z^\ell}{1-z^\ell}-z^{\ell
j}\right)=\s(z)-\frac{z^j}{1-z^j}.
\label{s_j}\end{equation}
\end{proposition}

\noindent{\bf Proof:} Let $a^{(k)}_j(n,m)$
denote the number  of Carlitz compositions of $n$ with the properties: 
\begin{itemize}
\item{} they have exactly $k$ parts
\item{} they do not contain a part of size $j$
\item{} the last part is of size $m$.
\end{itemize}
Then for $k\ge 1$ and $m\ne j$ we have
$$a^{(k+1)}_j(n,m)=a^{(k)}_j(n-m)-a^{(k)}_j(n-m,m),$$
where $a^{(k)}_j(\ell)$ is the number of Carlitz compositions of $\ell$
into $k$ parts, none of them equal to $j$.
Since $j$ is fixed throughout the argument, for the ease of notation
 we supress the subscript  $j$ throughout the argument. 
We additionally require that   $a^{(k)}(\ell)$ and $a^{(k)}(\ell,m)$
vanish whenever 
$\ell\le 0$. Let
$$f_k(z,u)=\sum_{m\ge1,\ n\ge0}a^{(k)}(n,m)z^nu^m.$$
Since the compositions enumerated in  $a^{(k)}(n,m)$ do not contain a part of size $j$, we have $a^{(k)}(n,j)=0$. Hence
\begin{eqnarray*}
f_{k+1}(z,u) & = & \sum_{n\ge0}\sum_{m\ge 1} a^{(k+1)}(n,m)z^nu^m \\
& = & \sum_{n\ge0}\sum_{m\ge 1\atop m\ne j} a^{(k+1)}(n,m)z^nu^m +\sum_{n\ge0}a^{(k+1)}(n,j)z^nu^j \\
& = & \sum_{n\ge0}\sum_{m\ge 1\atop m\ne j}\left(a^{(k)}(n-m)-a^{(k)}(n-m,m)\right)z^nu^m \\
& = & \sum_{n\ge0}\sum_{m\ge 1\atop m\ne j}a^{(k)}(n-m)z^nu^m 
\\ & & -
\sum_{n\ge0}\sum_{m\ge 1\atop m\ne j}a^{(k)}(n-m,m)z^nu^m 
\label{f(z,u)}
\end{eqnarray*}
The first sum above is equal to 
\begin{eqnarray*}& &\sum_{n\ge0,\
m\ge1}a^{(k)}(n-m)z^nu^m-\sum_{n\ge0}a^{(k)}(n-j)z^nu^j
\\ &=&\sum_{\ell\ge0,\
m\ge 1}a^{(k)}(\ell)(zu)^mz^\ell -\sum_{m\ge1} a^{(k)}(m)z^m(zu)^j \\ &=&
f_k(z,1)\frac{zu}{1-zu}-f_k(z,1)(zu)^j.
\end{eqnarray*}
where, in the second step we changed  the summation indices $n-m=\ell$
and $n-j=m$, respectively. By the same argument 
$$\sum_{m\ge1,\ \ell\ge0}a^{(k)}(\ell,m)z^\ell(zu)^m=f_k(z,zu).$$
Thus, for $k\ge 1$ we have
$$f_{k+1}(z,u)=f_k(z,1)\left(\frac{zu}{1-zu}-(zu)^j\right)-f_k(z,zu).$$
Letting $f_0(z,u)=1$ the last line can be rewritten as 
$$f_{k+1}(z,u)=f_k(z,1)\left(\frac{zu}{1-zu}-(zu)^j\right)-f_k(z,zu)+\delta_{k,0},$$
for $k\ge 0$. 
Introducing the function 
$$F(z,u)=\sum_{k\ge1}f_k(z,u)$$
and summing over $k\ge0$ we obtain
$$F(z,u)=F(z,1)\left(\frac{zu}{1-zu}-(zu)^j\right)+\left(\frac{zu}{1-zu}-(zu)^j\right)-F(z,zu).$$
This equation can be iterated 
to yield
\begin{eqnarray*}
F(z,u)&=&F(z,1)\left(\sum_{\ell\ge1}(-1)^{\ell-1}\left(\frac{z^\ell
u}{1-z^\ell u}-(z^\ell
u)^j\right)\right) \\ 
& & +\sum_{\ell\ge1}(-1)^\ell\left(\frac{z^\ell
u}{1-z^\ell u}-(z^\ell u)^j\right).
\end{eqnarray*}
Hence, for $u=1$ we get 
$$F(z,1)=F(z,1)\s_j(z)+\s_j(z),$$
where 
$\s_j(z)$ is defined by 
(\ref{s_j}).
Finally, letting 
$$C_j(z)=1+F(z,1),$$
we see that 
$$C_j(z)-1=\left(C_j(z)-1\right)\s_j(z)+\s_j(z),$$
which, since
$$F(z,1)=\sum_{k\ge1}f_k(z,1)=\sum_{k\ge1}\sum_{n\ge0\atop
m\ge1}a^{(k)}(n,m)z^n=\sum_{n\ge 0}a_{n,j}z^n,$$
 means that 
$$C_j(z)=\frac1{1-\s_j(z)},$$
is the generating function of the sequence $(a_{n,j})$.

\section{Singularities of the generating function}

A starting point of our analysis is the fact that the generating 
function
 of Carlitz compositions has the unique singularity in the
disc $\{z:\ |z|\le 0.663\}$. This singularity is the unique real root, $\rho$,
of the equation  $\s(z)=1$  on $[0,1]$. The numerical approximation of that root
is $\rho=0.571349\dots$.   Since $\s(z)$ and $\s_j(z)=\s(z)-z^j/(1-z^j)$, do
not differ by too much, 
functions $\s_j$ will have the same feature, at least for
$j$'s sufficiently large. In fact,   on the disc $\{z:\ |z|\le 0.663\}$ 
 the functions satisfy the following: 
there exists  $\delta>0$  such that
for every $j\ge 6$ the equation $\s_j(z)=1$ has the unique real
simple root on $[0,1]$.  Furthermore, these roots, which  will be denoted   by
$\rho_j$,   
have the following properties 
\begin{enumerate}
\item{} $\forall j\ge 6$, $0<\rho_j\le \rho+\delta$,
\item{} $\forall j\ge 6$ all  roots $\xi$ of
$\s_j(z)=1$  other than $\rho_j$ 
satisfy $|\xi|\ge \rho+2\delta$.
\item{} $\rho_j$ are strictly decreasing for $j\ge6$ and $\rho_j\to \rho$ as $j\to\infty$.
\end{enumerate}
A justification as well as a discussion of the remaining case $1\le j\le5$ is postponed until the appendix. 
We will need  asymptotics of $\rho_j$ and we will use the \lq\lq
bootstrapping method\rq\rq. Rewriting 
$$\s_j(z)=1$$ as 
\begin{equation}\s(z)=\frac1{1-z^j},
\label{rootj}\end{equation}
letting $\rho_j=\rho+\eps_j$, where $\eps_j=o(1)$ 
and substituting the latter expression for
$\rho_j$ into (\ref{rootj}) we get 
$$\s(\rho+\eps_j)=\frac1{1-(\rho+\eps_j)^j}=\frac1{1-\rho^j\left(1+\eps_j/\rho
\right)^j}.$$ 
Using Taylor's expansion on both sides we see that 
$$\s(\rho)+\s'(\rho)\eps_j+O(\eps_j^2)=
1+\rho^j\left(1+\frac{\eps_j}{\rho}\right)^j+O(\rho^{2j}).$$ 
Since $\s(\rho)=1$ we obtain 
that $\eps_j=\rho^j/\s'(\rho)+o(\rho^j)$. Hence 
\begin{equation}
\rho_j=\rho+\frac{\rho^j}{\s'(\rho)}+o(\rho^j),
\label{rho_j}
\end{equation} which is
sufficient for our purpose.
Let $A_j=-1/\s'_j(\rho_j)$ be the residue of $1/(1-\s_j(z))$ at
$\rho_j$. Then
$$\frac1{1-\s_j(z)}-\frac{A_j}{z-\rho_j}$$
is analytic in the disc $\{z:\ |z|\le \rho+\delta\}$ and by the Cauchy
integral formula we get
$$a_{n,j}= -\frac{A_j}{\rho_j}(\rho_j)^{-n}+O\left(\frac1{(\rho+\delta)^n}\right).$$
Finally, since $$\s_j(z)=\s(z)-\frac{z^j}{1-z^j},$$
we get 
$$\s_j'(\rho_j)=\s'(\rho_j)-\frac{j\rho_j^{j-1}}{(1-\rho_j^j)^2}.$$
Hence, taking into account (\ref{rho_j}) we get 
$$A_j=-\frac1{\s'(\rho)}+O(j\rho^j).$$
Consequently, for $j\ge 6$
\begin{equation}a_{n,j}=\left(\frac1{\s'(\rho)}+O(j\rho^j)\right)\frac1{\left(\rho+\rho^j\frac{1+o(1)}{\s'(\rho)}\right)^{n+1}}+O((\rho+\delta)^{-n}),
\label{anj}
\end{equation}
for some $\delta>0$ (universal for $j\ge 6$).

\section{Asymptotics}
The  following claim will account for most of the asymptotic analysis of (\ref{edn})
\begin{lemma} As $n\to\infty$,
$$\sum_{j=1}^{n}\left(1-\frac{a_{n,j}}{a_n}\right)=\sum_{j=0}^\infty\left(1-\left(1-\rho^j\a\right)^n\right)+o(1).\label{sum}$$
\end{lemma}

This statement immediately implies our theorem since the asymptotic behavior of the series on the right is known (see e.g.  \cite{FGD}, \cite{HitSte}, \cite[Section 7.8 and references therein]{SedFla}).

\noindent {\bf Proof of Lemma 3:}
The sequence $(a_n)$ was studied in \cite{KP} and one has
$$a_n=\frac1{\s'(\rho)}\left(\frac1{\rho}\right)^{n+1}+O((\rho+\delta)^{-n}).$$
Combining this with (\ref{anj}), for $j\ge 6$
$$\frac{a_{n,j}}{a_n}=
\frac{1+O(j\rho^j)}{\left(1+\rho^{j-1}\frac{1+o(1)}{\s'(\rho)}\right)^{n+1}}+
O((1+\Delta)^{-n}),
$$
where $\Delta=\delta/\rho>0$.
Consequently,
\begin{eqnarray*}
\sum_{j=1}^n\left(1-\frac{a_{n,j}}{a_n}\right)
&\!\!\!\!\!=\!\!\!\!\!&\sum_{j=1}^5\left(1-\frac{a_{n,j}}{a_n}\right) \\ 
&\!\!\!\!\!\!\!\!\!\!&\ \  +\sum_{j=6}^n\left(1-\frac{1+O(j\rho^j)}{\left(1+\rho^{j-1}\frac{1+o(1)}{\s'(\rho)}\right)^{n+1}}+O((1+\Delta)^{-n})\right).
\end{eqnarray*}
We will show that
$$\sum_{j=6}^n\frac{O(j\rho^j)}{(1+\rho^{j-1}\a_j)^{n+1}}=o(1),$$
where we have set $\alpha_j=1/\s'(\rho)+o(1)$. To this end, we will split this sum as 
$$\left(\sum_{j=6}^{j_n}+\sum_{j=j_n+1}^n\right)\frac{O(j\rho^j)}{(1+\rho^{j-1}\a_j)^{n+1}},$$
where $j_n$ will be chosen momentarily. Since $\a_j\to1/\s'(\rho)>0$ we have
$$\sum_{j=j_n+1}^n\frac{O(j\rho^j)}{(1+\rho^{j-1}\a_j)^{n+1}}\le C\sum_{j>j_n}j\rho^j=O(j_n\rho^{j_n})=o(1),$$
as long as $j_n\to\infty$ at {\em any} rate. The other term of the sum is bounded by
$$C\sum_{j=6}^{j_n}(1+\rho^{j-1}\a_j)^{-n}\le Cj_n\exp(-n\ln(1+\rho^{j_n-1}\a_j)).$$
Since $\ln(1+x)\ge x-x^2/2$, for $x\ge0$ we get the bound
$$Cj_n\exp\left(-n\a\rho^{j_n-1}\left(1-\frac{\a\rho^{j_n-1}}2\right)\right)
\le Cj_n\exp(-cn\rho^{j_n-1}).$$ Choosing $j_n\sim\mu\log_{1/\rho}n$ we see that this expression is bounded by 
$$O(\log n\cdot e^{-cn^{1-\mu}})=o(1),$$ whenever $\mu<1$.

It remains to consider the sum 
$$\sum_{j=5}^{n-1}\left(1-\left(1-\frac{\rho^j\alpha_j}{1+\rho^j\alpha_j}\right)^n\right).$$
We first replace this sum by a more convenient one
$$\sum_{j=5}^{n-1}\left(1-(1-\rho^j\alpha_j)^n\right).$$
The difference is at most 
$$\sum_{j=5}^\infty\left(\left(1-\frac{\rho^j\alpha_j}{1+\rho^j\alpha_j}\right)^n-(1-\rho^j\alpha_j)^n\right).$$
For a $j_0$ which will be chosen momentarily, we consider 
$$\sum_{j=5}^{j_0}\left(\left(1-\frac{\rho^j\alpha_j}{1+\rho^j\alpha_j}\right)^n-(1-\rho^j\alpha_j)^n\right),$$Since
each summand is nonnegative, 
this sum can be upperbounded by 
\begin{eqnarray*}\sum_{j=5}^{j_0}\left(1-\frac{\rho^j\alpha_j}{1+\rho^j\alpha_j}\right)^n
&\le&\sum_{j=5}^{j_0}\exp\left(-\frac{\alpha_j\rho^jn}{1+\rho^j\alpha_j}\right)
\\ &\le&
(j_0+1)\exp\left(-\frac{\alpha\rho^{j_0}n}{1+\rho^{j_0}\alpha_j}\right)\le(j_0+1)\exp\left(-c\rho^{j_0}n\right),
\end{eqnarray*}
for some absolute constant $c$. 
By choosing $j_0$ so that 
$$c\rho^{j_0}n\ge\k\ln n,$$ i.e.
$$j_0\le \frac{\ln((cn)/(\k\ln n))}{\ln(1/\rho)},$$
we see that the sum up to $j_0$ is bounded by $c\log n/n^\k$.
 For the remaining range of $j$'s we write
\begin{eqnarray*} \left(1-\frac{\rho^j\alpha_j}{1+\rho^j\alpha_j}\right)^n
-(1-\rho^j\alpha_j)^n 
&\!\!\!\!\!=\!\!\!\!\!& 
\left(1-\rho^j\alpha_j+\frac{\rho^{2j}\alpha^2_j}{1+\rho^j\alpha_j}\right)^n
-(1-\rho^j\alpha)^n  \\
&\!\!\!\!\!\le \!\!\!\!\!&n\frac{\alpha^2\rho^{2j}}{1+\rho^j\alpha}
\left(1-\frac{\rho^j\alpha}{1+\rho^j\alpha}\right)^{n-1}\le 
n\alpha^2\rho^{2j},
\end{eqnarray*}
where we have used the inequality $(a+b)^n-a^n\le nb(a+b)^{n-1}$ valid for 
nonnegative numbers $a$ and $b$.
Hence 
$$\sum_{j\ge j_0}\left(\left(1-\frac{\rho^j\alpha_j}{1+\rho^j\alpha_j}\right)^n-
(1-\rho^j\alpha_j)^n\right)\le n\alpha^2_j\sum_{j\ge 
j_0}\rho^{2j}\le cn\rho^{2j_0}.$$
The choice of $j_0$ so that $\rho^{j_0}=\Theta(\log n/n)$ is within the 
previous constraint, and for that choice we have 
$$cn\rho^{2j_0}=\Theta(\frac{\log^2n}{n}).$$  
Using the same argument we can show that $\a_j$ can be replaced by its limit $\a=1/\s'(\rho)$, i.e. that we have 
$$\sum_{j=5}^{n-1}\left(1-(1-\rho^j\alpha_j)^n\right)=
\sum_{j=5}^{n-1}\left(1-(1-\rho^j\a)^n\right)+o(1).$$ 
Finally, the sum on the right can be increased 
to 
$$\sum_{j=0}^\infty\left(1-(1-\rho^j\a)^n\right),$$
since
$$\sum_{j\ge n}\left(1-(1-\rho^j\a)^n\right)\le\a n
\sum_{j\ge n}\rho^j
=o(1)$$
and is thus negligible. Also, for each fixed $j$, with $1\le j\le 5$,
 $a_{n,j}/a_n=o(1)$   as we will indicate in the appendix. This means that 
$$\sum_{j=1}^5\left(1-\frac{a_{n,j}}{a_n}\right)=5+o(1),$$
and proves the  Lemma.

\section{Appendix}

In order to show that $\s_j(z)=1$ has a unique real root $\rho_j$ and
the existence of a $\delta>0$ satisfying the asserted properties we
rewrite $\s(z)$ in a more convenient form
$$\s(z)=\sum_{m=1}^\infty\frac{z^m}{1+z^m},$$
which can be done by expanding $z^j/(1-z^j)$ into geometric series and
interchanging the order of summation. Now, $\s_j(z)=1$ can be
rewritten as 
$$\sum_{m=1}^\infty\frac{z^m}{1+z^m}-\frac{z^j}{1-z^j}-1=0.$$
We want to use Rouch\'e's theorem. To this end split the left hand side as $f(z)+g_j(z)$ where 
$$f(z)=\sum_{m=1}^6\frac{z^m}{1+z^m}-1.$$
It can be verified that 
$$\min_{|z|=0.663}|f(z)|\ge 0.28.$$
Actually Maple suggests that a stronger claim is true, namely that
minimum of $|f(z)|$ on that circle is attained at $z=0.663$ and is
$0.283467$.
We have not tried to prove it and thus we claim only the weaker
statement, which can be verified by evaluating $|f|$ at sufficiently
many points and using its Lipshitz property: for $|z_1|=|z_2|=r$,
\begin{eqnarray*}\Big|\ |f(z_1)|-|f(z_2)|\ \Big| &\le&|f(z_1)-f(z_2)|\le
\sum_{m=1}^6\frac{|z_1^m-z_2^m|}{|(1+z_1^m)(1+z_2^m)|} \\ 
&\le&
|z_1-z_2|\sum_{m=1}^6\frac{mr^{m-1}}{(1-r^m)^2}\le C|z_1-z_2|,
\end{eqnarray*}
where 
$C\le (1-r)^{-2}\sum_{m\ge 1}mr^{m-1}\le(1-r)^{-4}$. 
Since for $|z|=0.663$ 
\begin{eqnarray*}|g_j(z)|&\le&
\sum_{m=7}^\infty\frac{|z|^m}{1-|z|^m}+\frac{|z|^j}{1-|z|^j} \\ &\le& 
\frac{0.663^7}{(1-0.663^7)(1-0.663)}+\frac{0.663^6}{1-0.663^6}\le
0.27,
\end{eqnarray*}
we see that 
$$\min_{|z|=0.663}|f(z)|>\max_{|z|=0.663}|g_j(z)|,$$
uniformly for all $j\ge 6$. Thus $f+g_j$ and $f$ have the same number of roots
in the disc $|z|\le 0.663$. But $f$ has a unique (and thus necessarily real) root on that
disc (this again may be verified by Maple after  converting it to a polynomial equation). The existence
of $\delta$ follows since the same argument can be repeated for a disc
with a radius just a bit smaller. As for the monotonicity, it suffices
to rewrite $\s_j(z)=1$ as $\s(z)=1/(1-z^j)$ and notice that $\s(z)$ is  increasing as a function of real variable
and that 
$$\frac1{1-z^{j_1}}<\frac1{1-z^{j_2}},$$ whenever $j_1>j_2$ and $0<z<1$ is
real.

The same argument (with changed parameters) can be used to verify that 
 for $j=2,3,4$ and $5$, on the disc $|z|<r$, $0<r<1$,
 $\s_j(z)=1$ has a unique real root $\rho_j>\rho$. 
The hardest case is $j=2$. We found that the splitting 
$$f_2(z)=\sum_{m=1}^{25}\frac{z^m}{1+z^m}-\frac1{1-z^2},\quad g(z)=\sum_{m=26}^\infty\frac{z^m}{1+z^m},$$ 
will do the job. 
On  the circle $|z|=0.8$ one has $|f_2(z)|\ge 0.06$, $|g(z)|\le 0.016$. The polynomial resulting from multiplying $f_2$ by the product of denominators is of degree $264$ (after cancellations). Maple (somehow reluctantly) shows that the   root closest to zero is real and is about $0.78397$ (there is another real root of about $0.927122$) and the next  closest to zero roots are complex conjugate and have absolute value around $0.81914$. For $j=3,4,5$ one can get away with letting $g(z)=\sum_{m=11}^\infty z^m/(1+z^m)$,  and choosing $r=0.75$.

The same method could presumably be used to force the argument for $j=1$. However,  here matters would be computationally worse. Furthermore, the equation $\s_1(z)=1$ does not have real solutions, and thus, its closest to zero root (if it exists) would have to come from a pair of complex conjugates. Luckily, for that case we can use a different argument based on the probabilistic approach used in \cite{HS} (see also \cite{HitSte}). Let us briefly sketch it. Consider the set $C_n$ of {\em all} $2^{n-1}$ compositions of an integer $n$ and let $\Q$ be the uniform
probability measure on $C_n$. Since the restriction of such a measure to any subset is again the uniform measure on that subset we can view the uniform measure $\P$ on the set of all Carlitz compositions as a conditional measure obtained by restricting $\Q$ to $\Omega_n$, i.e. 
$$\P(\ \cdot\ )=\Q(\ \cdot\ |\Omega_n).$$
Let $A_1^c$ be the set of {\em all} compositions of $n$ that do not
use part size 1 and recall that $I_1^c$ is the set of all {\em
Carlitz} compositions with this property. Then,
$$\frac{a_{n,1}}{a_n}=\P(I_1^c) = \Q(A_1^c|\Omega_n)=
\frac{\Q(A_1^c\cap\Omega_n)}
{\Q(\Omega_n)} \le \frac{\Q(A_1^c)}{\Q(\Omega_n)}.
$$
Now, by the result of Knopfmacher and Prodinger 
$$\Q(\Omega_n)\ge c\frac{1.75^n}{2^n}=c(0.875)^n,$$
and we need to upper bound $\Q(A_1^c)$. To this end we will use the
observation made in \cite{HS} or \cite{HitSte} that a random
composition of $n$ is distributed like 
$$(\tG_1,\dots,\tG_\tau),$$ where, with  $(\G_i)$ being a sequence of i.i.d. geometric random variable with parameter $1/2$, we let 
$$\tau=\inf\{k:\ \sum_{q=1}^k\G_q\ge n\},$$ and 
$$\tG_q=\left\{ \begin{array}{ll}\G_q & 
\mbox{if $q<\tau$,} \\ & \\ 
n-\sum_{p=1}^{\tau-1}\G_q& \mbox{if $q=\tau$.}\end{array}\right.$$ 
Furthermore, $\tau$ is distributed like 1 plus  a binomial random variable with parameters $n-1$ and $1/2$. Therefore,
\begin{eqnarray*}
\Q(A_1^c)&=&\Q(\bigcap_{k=1}^\tau\{\tG_k>1\}) \le 
\Q(\bigcap_{k=1}^{\tau-1}\{\G_k>1\}) \\ 
&\le& \Q(\tau\le k_0)+\Q(\bigcap_{k=1}^{\tau-1}\{\G_k>1\}\cap\{\tau>k_0\}) \\ &
\le& \Q(\tau -\E\tau\le k_0-\E\tau)+\Q(\bigcap_{k=1}^{k_0}\{\G_k>1\}) \\ &\le&
\Q(|\tau-\E\tau|>\frac{n+1}2-k_0)+\left(1-\frac12\right)^{k_0} \\ &\le& 
2\exp\left(-\frac{((n+1)/2-k_0)^2}{2(n-1)/4}\right)+\exp\left(-k_0\ln2\right) \\ &\le& 2\exp\left(-2\frac
{\phi_n^2}{n-1}\right)+\exp\left(-k_0\ln2\right),
\end{eqnarray*}
where we have put $k_0=(n+1)/2-\phi_n$ and $\phi_n=\alpha(n+1)$, for some $0<\alpha<1/2$.
Then we get 
$$
\Q(A_1^c)\le C\left(e^{-2\alpha^2}\right)^n+C\left(e^{-(1/2-\alpha)\ln2}\right)^n.$$
To be able to claim that $\Q(A_1^c)/\Q(\Omega_n)$ tends to zero, we need both of the inequalities 
$$e^{-2\alpha^2}<0.875\quad\hbox{and}\quad e^{-(1/2-\alpha)\ln2}<0.875
$$ 
to be satisfied simultaneously. But this can be achieved by choosing $\alpha$
to be any number subject to $0.26<\alpha<0.3$. This completes the
argument. 

\vskip .5cm
\noindent
{\bf Acknowledgment} This papers owes its existence to the curiosity of Herbert Wilf, who
brought the paper \cite{KP} to our attention by 
asking the second named author whether the methods used in \cite{HS}
could be used to rederive some of the results of \cite{KP} (that question, by
the way, 
remains unanswered).

\end{document}